\documentclass[oneside,10pt]{article}
\usepackage[b5paper]{geometry}

\usepackage{amsmath,amsfonts,amsxtra,amssymb,latexsym,amscd,amsthm}

\DeclareMathOperator{\tr}{tr}
\DeclareMathOperator{\rank}{rank}
\DeclareMathOperator{\Ima}{Im}
\DeclareMathOperator{\Rea}{Re}

\DeclareMathOperator{\per}{per}
\DeclareMathOperator{\Ker}{Ker}

\usepackage{tikz}
\usepackage{array}
\usepackage{dsfont}
\usepackage{enumitem}
\usepackage{amsthm} 
\usepackage{mathrsfs,upref}       
\usepackage{mathptmx}

\newtheorem{theorem}{Theorem}           
\newtheorem{lemma}{Lemma}               

\newtheorem{proposition}{Proposition}[section]

\theoremstyle{definition}

\newtheorem{definition}{Definition}

\newtheorem{remark}{Remark}
\newtheorem{Conjecture}{Conjecture}
\numberwithin{equation}{section} 

\usepackage{mathrsfs}

\begin{document}
\title{A simple counterexample for the permanent-on-top conjecture}
\author{Tran Hoang Anh\\Eötvös Loránd University. Institute of Mathematics.\\E-mail: bongtran5399@gmail.com\\ }
\maketitle

\begin{abstract}
The permanent-on-top conjecture (POT) was an important conjecture on the largest eigenvalue of the Schur power matrix of a positive semi-definite Hermitian matrix, formulated by Soules. The conjecture claimed that for any positive semi-definite Hermitian matrix $H$, $\per(H)$ is the largest eigenvalue of the Schur power matrix of the matrix $H$. After half a century, the POT conjecture has been proven false by the existence of  counterexamples which are checked with the help of computer. It raises concerns about a counterexample that can be checked by hand (without the need of computers). A new simple counterexample for the permanent-on-top conjecture is presented which is a complex matrix of dimension 5 and rank 2.

\end{abstract}

\section{Introduction and notations}
The symbol $S_n$ denotes the symmetric group on $n$ objects. The permanent of a square matrix is a vital function in linear algebra that is similar to the determinant. For an $n \times n$ matrix $A=(a_{ij})$ with complex coefficients, its permanent is defined as $\per(A)=\sum_{\sigma \in S_n}\prod_{i=1}^{n}a_{i,\sigma(i)}$. By $\mathscr{H}_n$ we mean the set of all $n \times n$ positive semi-definite Hermitian matrices. The Schur power matrix of a given $n \times n$ matrix $A=(a_{ij})$, denoted by $\pi(A)$, is a $n! \times n!$ matrix with the elements indexed by permutations $\sigma, \tau \in S_n$: $$\pi_{\sigma \tau}(A)= \prod_{i=1}^{n}a_{\sigma(i) \tau(i)}$$
\begin{Conjecture}$\textbf{The permanent-on-top conjecture (POT) \cite{Soules}:}$
	Let $H$ be an $n \times n$ positive semi-definite Hermitian matrix, then $\per(H)$ is the largest eigenvalue of $\pi(H)$.
\end{Conjecture}
In 2016, Shchesnovich provided a 5-square, rank 2 counterexample to the permanent-on-top conjecture with the help of computer \cite{Shchesnovich}.
\begin{definition}
	For an $n \times n$ matrix $A=(a_{ij})$, let $d_A$ be a function $S_n \rightarrow \mathbb{C}$ defined by $$d_A(\sigma)=\prod_{i=1}^{n}a_{\sigma(i)i}$$
	This function is also called the "diagonal product" function \cite{BS}. Then we  can define $\det(A)=\sum_{\sigma \in S_n}(-1)^{sign(\sigma)}d_A(\sigma)$ and $\per(A)= \sum_{\sigma \in S_n}d_A(\sigma)$. 
\end{definition}
For any $n$-square matrix $A$ and $I, J \subset [n]$, $A[I,J]$ denotes the submatrix of $A$ consisting of entries which are the intersections of $i$-th rows and $j$-th columns where $i \in I,\ j \in J$. We define $A(I,J)=A[I^c,J^c]$.\\
In this paper, we shall study the properties of the spectrum of the Schur power matrix by examining the spectra of the matrices $\mathscr{C}_k(A)$ which are defined in the manner:\\
	For any $1\leq k \leq n$, the matrix $\mathscr{C}_k(A)$ is a matrix of size $ \binom{n}{k} \times \binom{n}{k}$ with its $(I,J)$ entry ($I$ and $J$ are $k$-element subsets of $[n]$) defined by $\per(A[I,J]).\per(A[I^c,J^c])$. There is another conjecture on these matrices $\mathscr{C}_k(A)$ which states that:
\begin{Conjecture}\textbf{Pate's conjecture \cite{PT}}
	Let $A$ be an $n \times n$ positive semi-definite Hermitian matrix and $k$ be a positive integer number less than $n$, then the largest eigenvalue of $\mathscr{C}_k(A)$ is $\per(A)$.	
\end{Conjecture}
Pate's conjecture is weaker than the permanent-on-top conjectu
POT) because it is well-known that every eigenvalue of $\mathscr{C}_k$ is also an eigenvalue of the Schur power matrix. In the case $k=1$, in \cite{BS}, it was conjectured that $\per(A)$ is necessarily the largest eigenvalue of $\mathscr{C}_1(A)$ if $A \in \mathscr{H}_n$. Stephen W. Drury has provided an 8-square matrix as a counterexample for this case in the paper \cite{DS1}. Besides, Bapat and Sunder raise a question as follows:
\begin{Conjecture}\textbf{Bapat \& Sunder conjecture:}
	Let $A$ and $B=(b_{ij})$ be $n \times n$ positive semi-definite Hermitian matrices, then 
	$$\per(A \circ B) \leq \per(A)\prod_{i=1}^{n}b_{ii}$$ where $A \circ B$ is the entrywise product (Hadamard product).
\end{Conjecture}
The Bapat \& Sunder conjecture is weaker than the permanent-on-top conjecture and has been proved false by a counterexample which is a positive semi-definite Hermitian matrix of order 7 proposed by Drury \cite{DS2}. In the present paper, a new simple counterexample for the permanent-on-top conjecture and Pate's conjecture is presented. It has size $5 \times 5$ and rank 2.
\begin{Conjecture}\textbf{The Lieb permanent dominance conjecture 1966 \cite{Lieb}}
		Let $H$ be a subgroup of the symmetric group $S_n$ and let $\chi$ be a character of degree $m$ of $H$. Then $$\dfrac{1}{m}\sum_{\sigma \in H}\chi(\sigma)\prod_{i=1}^{n}a_{i\sigma(i)} \leq \per(A)$$ holds for all $n \times n$ positive semi-definite Hermitian matrix A.
\end{Conjecture}
The permanent dominance conjecture is weaker than the permanent-on-top conjecture and still open. The POT conjecture was proposed by Soules in 1966 as a strategy to prove the permanent dominance conjecture.
\begin{definition}
The elementary symmetric polynomials in $n$ variables $x_1,x_2,\dotsc,x_n$ are $e_k$ for $k=0,1,\dotsc,n$. In this paper, we define $e_k(x_i)$ for $i=1,2,\dotsc,n$ to be the elementary symmetric polynomial of degree $k$ in $n-1$ variables obtained by erasing variable $x_i$ from the set $\lbrace x_1,x_2,\dotsc,x_n\rbrace$ and, for any sebset $I\subset [n]$, the notation $e_k[I]$ denote the elementary symmetric polynomial of degree $k$ in $|I|$ variables $x_i$'s, $i \in I$.\\

\end{definition}
\section{Associated matrices}
We define the associated matrix of a matrix representation $W: S_n \rightarrow GL_N(\mathbb{C})$ with respect to a $n \times n$ matrix $A$ by: $$M_{W}(A)=\sum_{\sigma \in S_n}d_A(\sigma)W(\sigma)$$
\begin{proposition}
	The Schur power matrix of a given $n \times n$ Hermitian matrix $A$ is the associated matrix of the left-regular representation with respect to $A$.
\end{proposition}
\begin{proof}
	Take a look at the $(\sigma,\tau)$ entry of $M_{L}(A)$ which is
	 $$\sum_{\eta \in S_n,\ \eta \circ \tau=\sigma}d_A(\eta)=d_A(\sigma \circ \tau^{-1})=\prod_{i=1}^{n}a_{\sigma(i)\tau(i) }$$
	 the right side is the $(\sigma,\tau)$ entry of $\pi(A)$.
\end{proof}
Let us now consider two important matrices $\mathscr{C}_1(A)$ amd $\mathscr{C}_2(A)$ that shall appear frequently from now on.

\begin{definition}
	Let $\mathscr{N}_k: S_n \rightarrow GL_{\binom{n}{k}}(\mathbb{C})$ be the matrix representation given by the permutation action of $S_n$ on $\binom{[n]}{k}$.
\end{definition}
\begin{proposition}
		For any $n \times n$ Hermitian matrix $A$, the matrix $\mathscr{C}_k(A)$ is the matrix $M_{\mathscr{N}_k}(A)$.
\end{proposition}
We obtain directly the statement that every eigenvalue of matrix $M_{\mathscr{N}_k}(A)$ is an eigenvalue of the associated matrix of the left-regular representation which is the Schur power matrix. Consequently, Pate's conjecture is weaker than the permanent-on-top conjecture(POT).\\
\section{Several properties of the Schur power matrix and $\mathscr{C}_1(A)$ in rank 2 case}
The main object of this section is $n \times n$ positive semi-definite Hermitian matrices of rank 2. We know that every matrix $A \in \mathscr{H}_n$ of rank 2 can be written as the sum $v_1v_1^*+v_2v_2^*$ where $v_1$ and $v_2$ are two column vectors of order $n$.
\begin{definition}
	A matrix $A \in \mathscr{H}_n$ is called  "formalizable" if $A$ can be written in the form $v_1v_1^*+v_2v_2^*$ and every element of $v_1$ vector is non-zero.
\end{definition}
\begin{definition}
	The formalized matrix $A'$ of a given formalizable matrix $A$ defined in the manner:
	if $A=v_1v_1^*+v_2v_2^*$ and $v_1=(a_1,\dotsc,a_n)^T, \ a_i \neq 0\  \forall i=1,\dotsc,n$ ; $v_2=(b_1,\dotsc,b_n)^T$ then $A'=v_3v_3^*+v_4v_4^*$ where $v_3=(1,\dotsc,1)^T$ and $v_4=(\frac{b_1}{a_1},\dotsc,\frac{b_n}{a_n})^T$.
	
\end{definition}
\begin{proposition}: Let $A \in \mathscr{H}_n$ be a formalizable matrix, then $\pi(A)=\prod_{i=1}^{n}|a_i|^2 \pi(A')$.
\end{proposition}
\begin{proof}
We compare the $(\sigma,\tau)$-th entries of two matrices.
$$\pi_{\sigma \tau}(A)=\prod_{i=1}^{n}(a_{\sigma(i)}\overline{a_{\tau(i)}}+b_{\sigma(i)}\overline{b_{\tau(i)}})=\prod_{i=1}^{n}|a_i|^2\prod_{i=1}^{n}\left(1+\frac{b_{\sigma(i)}}{a_{\sigma(i)}}\frac{\overline{b_{\tau(i)}}}{\overline{a_{\tau(i)}}}\right)=\prod_{i=1}^{n}|a_i|^2\pi_{\sigma \tau}(A')$$
\end{proof}
\begin{remark}
	The same result will be obtained with the matrices $\mathscr{C}_k(A)$ and $\mathscr{C}_k(A')$. It is obvious to see that if the matrix $A$ is a counterexample for the permanent-on-top conjecture and Pate's conjecture then so is $A'$. Assume that we have an unformalizable matrix $B \in \mathscr{H}$ of rank 2 that is a counterexample for the permanent-on-top conjecture and Pate's conjecture. That also implies that there is a column vector $x$ such that the following inequality holds $$\dfrac{x^*\pi(B)x}{\|x\|^2}>\per(B)$$
	By continuity and $B=vv^*+uu^*$, we can change slightly the zero elements of the vector $v$ such the the inequality remains. Therefore, if the permanent-on-top conjecture or Pate's conjecture is false for some positive semi-definite Hermitian matrix of rank 2 then so is the permanent-on-top conjecture and Pate's conjecture for some formalizable matrices. That draws our attention to the set of all formalizable matrices.\\
	For any $n \times n$ positive semi-definite Hermitian matrix $A$ of rank 2 there exist two eigenvectors of $v$ and $u$ of $A$ such that $A=vv^*+uu^*$. Let $u_i,v_i$ be the $i$-th row elements of $v$ and $u$ respectively for $i=\overline{1,n}$. In the case $A$ has a zero row then $\per(A)=0$ and the Schur power matrix and matrices $\mathscr{C}_k(A)$ of $A$ are all zero matrices, there is nothing to discuss. Otherwise, every row of $A$ has a non-zero element (so does every column since $A$ is a Hermitian matrix) which means that for any $i=\overline{1,n}$, the inequalities $|v_i|^2+|u_i|^2>0$ hold. Besides, $A$ can be rewritten in the form $$(\sin(x)v+\cos(x)u)(\sin(x)v+\cos(x)u)^*+(\cos(x)v-\sin(x)u)(\cos(x)v-\sin(x)u)^* \ \forall x \in [0,2\pi]$$ and the system of $n$ equations $\sin(x)v_i+\cos(x)u_i=0,\ i=\overline{1,n}$ takes finite solutions in the interval $[0,2\pi]$. Therefore, there exists $x \in [0,2\pi]$ satisfying that $(\sin(x)v+\cos(x)u)$ has every element different from 0. Hence, every rank 2 positive semi-definite Hermitian that has no zero-row is formalizable. Several properties about the formalized matrices are presented below.
\end{remark}
Let $H \in \mathscr{H}_n$  be a formalizable matrix of the form $H=vv^*+uu^*$ where $v=(1,\dotsc,1)^T$ and $u=(x_1,x_2,\dotsc,x_n)^T$. We recall quickly the Kronecker product \cite{Zhang1}.
\begin{definition}
	The Kronecker product (also known as tensor product or direct product) of two matrices $A$ and $B$ of sizes $m \times n$ and $s \times t$, respectively, is defined to be the $(ms)\times (nt)$ matrix
	$$A \otimes B=
	\begin{pmatrix}
	a_{11}B&a_{12}B&\dotsc a_{1n}B\\
	a_{21}B&a_{22}B&\dotsc a_{2n}B\\
	\vdots&\vdots&\vdots\\
	a_{n1}B&a_{n2}B&\dotsc a_{nn}B\\
	\end{pmatrix}$$
\end{definition}
\begin{lemma}\textbf{The upper bound of rank of the Schur power matrix of rank 2: }
	If $A$ is  $n \times n$ of rank 2 then rank of $\pi(A)$ is not larger than $2^n-n$.	
\end{lemma}
\begin{proof}
	We observe that $\rank(A)=2$ implies that $\dim(\Ima(A))=2$ and $\dim(\Ker(A))=n-2$. Let $\langle w,t \rangle$ be an orthonormal basis of the orthogonal complement of $\Ker(A)$ in $\mathbb{C}^n$, then denote $v=Aw, \ u=At$. Thus, $A$ can be rewritten in the form $vw^*+ut^*$ where $v=(a_1,\dotsc,a_n)^T,\ u=(b_1,\dotsc,b_n)^T$. It is obvious that $\Ima(A)=\langle v,u \rangle$. Let us denote the Kronecker product of $n$ copies of the matrix $A$ by $\otimes^nA$. The mixed-product property of Kronecker product implies that $\Ima(\otimes^nA)=\langle\{\otimes_{i=1}^{n}t_i,\ t_i \in \{v,u\}\}\rangle$. Furthermore, the Schur power matrix of $A$ is a diagonal submatrix of $\otimes^nA$ obtained by deleting all entries of $\otimes^nA$ that are products of entries of $A$ having two entries in the same row or column. Let define the function $f$ in the manner that $$f: \{\otimes_{i=1}^{n}t_i,\ t_i \in \{v,u\}\} \rightarrow \widetilde{V}$$
	and the $\sigma$-th element of $f(\otimes_{i=1}^{n}\alpha_i)$ vector of order $n!$ is $\prod_{i=1}^{n}t_i(\sigma(i))$ where $t_i(j)$ is the $j$-th row element of the column vector $t_i$. Let $\mathscr{B}=\{f(\otimes_{i=1}^{n}t_i),\ t_i \in \{v,u\} \}$ then $\mathscr{B}$ is a generator of $\Ima(\pi(A))$ since $\pi(A)$ is a principal matrix of $\otimes^nA$ and $\Ima(\otimes^nA)=\langle\{\otimes_{i=1}^{n}t_i,\ t_i \in \{v,u\}\}\rangle$. We partition $\mathscr{B}$ into disjoint sets $S_k$ $$k=0,1,\dotsc,n, \ S_k=\{f(\otimes_{i=1}^{n}t_i),\ t_i \in \{v,u\}, \ \text{v appears k times in the Kronecker product}\}$$
	Hence, for any $k=1,2,\dotsc,n$ the $\sigma$-th row element of the sum vector $\sum_{w \in S_k}w$ is $$\sum_{\substack{1\leq i_1<\dotsc<i_k\leq n \\ 1\leq i_{k+1}<\dotsc<i_n\leq n}}\prod_{j=1}^{k}a_{\sigma(i_j)}\prod_{t=k+1}^{n}b_{\sigma(i_t)}
	=\sum_{\substack{1\leq i_1<\dotsc<i_k\leq n \\ 1\leq i_{k+1}<\dotsc<i_n\leq n}}\prod_{j=1}^{k}a_{i_j}\prod_{t=k+1}^{n}b_{i_t}$$ and $S_0=\{(1,,1,\dotsc,1)^T\}$. Therefore, for any $k=1,\dotsc,n$ then $S_0 \cup S_k$ is linearly dependent. Hence, by deleting an arbitrary element of each set $S_k$ $,\ k=1,\dotsc,n$, then it still remains a generator of $\Ima(\pi(H))$. Thus
	$$\rank(\pi(A))=\dim(\Ima(\pi(A)))\leq |\mathscr{B}|-n=2^n-n$$
	
\end{proof}
\begin{lemma}\textbf{The permanent of a formalized matrix \cite{Marcus}:}
	$$\per(H)=\sum_{k=0}^{n}k!(n-k)!|e_k|^2$$
\end{lemma}
\begin{proof}
	We show that $$\begin{aligned}
		\per(H)&=\sum_{\sigma \in S_n}\prod_{i=1}^{n}(1+x_{i}\overline{x_{\sigma(i)}})\\
		&=n!+\sum_{\sigma \in S_n}\sum_{k=1}^{n}\sum_{1\leq i_1<\dotsc< i_k\leq n}x_{i_1}\dotsc x_{i_k}\overline{x_{\sigma(i_1)}\dotsc x_{\sigma(i_k)}}\\
		&=n!+\sum_{k=1}^{n}\sum_{1\leq i_1<\dotsc< i_k\leq n}x_{i_1}\dotsc x_{i_k}\overline{\sum_{\sigma \in S_n}x_{\sigma(i_1)}\dotsc x_{\sigma(i_k)}}\\
		&=n!+\sum_{k=1}^{n}\sum_{1\leq i_1<\dotsc< i_k\leq n}k!(n-k)!x_{i_1}\dotsc x_{i_l}\overline{e_k}\\
		&=\sum_{k=0}^{n}k!(n-k)!|e_k|^2\\
	\end{aligned}$$
\end{proof}
We use the elementary symmetric polynomials to examine entries of $\mathscr{C}_1(H)$ with the $(i,j)$-th entry defined by $(1+x_i\overline{x_j}).\per(H(i|j))$ and
$$\begin{aligned}
	\per(H(i|j))&=\sum_{\sigma \in S_n;\ \sigma(i)=j}\prod_{l\neq i}(1+x_l\overline{x_{\sigma(l)}})\\
	&=\sum_{\sigma \in S_n;\ \sigma(i)=j}\sum_{k=0}^{n-1}\sum_{1\leq i_1< \dotsc < i_k \leq n;\ i_m\neq i \ \forall m=1,\dotsc,k}x_{i_1}\dotsc x_{i_k}\overline{x_{\sigma(i_1)}\dotsc x_{\sigma(i_k)}}\\
	&=\sum_{k=0}^{n}\sum_{1\leq i_1<\dotsc< i_k\leq n,\ i_m\neq i}k!(n-1-k)!x_{i_1}\dotsc x_{i_l}\overline{e_k(x_j)}\\
	&=\sum_{k=0}^{n-1}k!(n-1-k)!e_k(x_i)\overline{e_k(x_j)}\\
\end{aligned}$$
And notice that $$e_k=x_ie_{k-1}(x_i)+e_k(x_i) \ \forall k=1,\dotsc,n$$
Then $$\begin{aligned}
	\dfrac{\per(H)}{n}&=\frac{1}{n}\sum_{k=0}^{n}k!(n-k)!|e_k|^2\\
	&=(n-1)!(|e_0|^2+|e_n|^2)+\sum_{k=1}^{n-1}\dfrac{k!(n-k)!}{n}(x_ie_{k-1}(x_i)+e_k(x_i))\overline{(x_je_{k-1}(x_j)+e_k(x_j))}\\
\end{aligned}$$
Hence $$\begin{aligned}
&\quad (1+x_i\overline{x_j}).\per(H(i|j))-\dfrac{\per(H)}{n}\\ &=\sum_{k=1}^{n-1}\left(k!(n-1-k)!-\frac{k!(n-k)!}{n}\right) e_k(x_i)\overline{e_k(x_j)}\notag\\
&\quad+\left((k-1)!(n-k)!-\frac{k!(n-k)!}{n}\right)x_ie_{k-1}(x_i)\overline{x_je_{k-1}(x_j)}\notag\\
&\quad-\frac{k!(n-k)!}{n}(x_ie_{k-1}(x_i)\overline{e_k(x_j)}+\overline{x_je_{k-1}(x_j)}e_k(x_i)) \\
&=\sum_{k=1}^{n-1}\frac{(k-1)!(n-1-k)!}{n}(ke_k(x_i)-(n-k)x_ie_{k-1}(x_i))(k\overline{e_k(x_j)}-(n-k)\overline{x_je_{k-1}(x_j)})\\
&=\sum_{k=1}^{n-1}\frac{(k-1)!(n-1-k)!}{n}(ne_k(x_i)-(n-k)e_k)(\overline{ne_k(x_j)-(n-k)e_k})
\end{aligned}$$
Therefore, we have the following proposition.
\begin{proposition}
 The matrix $\mathscr{C}_1(H)$ can be rewritten in the form
$$\mathscr{C}_1(H)=\frac{\per(H)}{n} vv^*+\sum_{k=1}^{n-1}\frac{(k-1)!(n-1-k)!}{n}v_kv_k^*$$
where $v=(1,\dotsc,1)^T$ of order $n$, for $k=1,\dotsc,n-1$, $v_k=(\dotsc,\underset{\text{i-th element}}{\underbrace{ne_k(x_i)-(n-k)e_k}},\dotsc)^T$
\end{proposition}
\begin{proposition}
	For any $k=1,\dotsc,n-1$, $\langle v,v_k\rangle=0$
\end{proposition}
\begin{proof}
	$$\begin{aligned}
		\langle v,v_k\rangle&=\sum_{i=1}^{n}(ne_k(x_i)-(n-k)e_k)\\
		&=n\sum_{i=1}^{n}e_k(x_i)-n(n-k)e_k\\
		&=0\\
		\end{aligned}$$
\end{proof}
\begin{proposition} The rank of $\mathscr{C}_1(H)$ is the cardinality of the set $\{x_i,i=\overline{1,n}\}$. In formula, $\rank(\mathscr{C}_1(H))=|\{x_i,i=\overline{1,n}\}|$.
\end{proposition}
\begin{proof}
	For the $i$-th element of $v_k$, we have 
	$$\begin{aligned}
		&\quad ne_k(x_i)-(n-k)e_k\\
		&=ke_k-nx_ie_{k-1}(x_i)\\
		&=ke_k+n\sum_{j=1}^{k}(-1)^je_{k-j}x_i^j\\
	\end{aligned}$$
which leads us to a conclusion that $\langle v,v_1,\dotsc,v_{n-1}\rangle=\langle p_0,\dotsc,p_{n-1}\rangle$ where $p_j=(\dotsc,\underset{\text{i-th element}}{\underbrace{x_i^j}},\dotsc)^T$ which is equal to $|\{x_i,i=\overline{1,n}\}|$ by the determinantal formula of Vandermonde matrices.
\end{proof}
\begin{proposition} The determinant of $\mathscr{C}_1(H)$ is given by
	$$\det(\mathscr{C}_1(H))=\frac{\per(H)}{n}\prod_{k=1}^{n-1}n(k-1)!(n-1-k)!\cdot \prod_{i<j}|x_i-x_j|^2$$

\end{proposition}
\begin{proof}	
	Case 1: There are indices $i$ and $j$ such that $x_i=x_j$ then $\rank(\mathscr{C}_1(H))<n$ that is equivalent to $\det(\mathscr{C}_1(H))=0$.\\
	Case 2: $x_i$'s are distinct then $\{v,v_1,\dotsc,v_{n-1}\}$ makes a basis of $\mathbb{C}^n$. Therefore, $\mathscr{C}_1(H)$ is similar to the Gramian matrix of $n$ vectors\\ $\left\{\sqrt{\frac{\per(H)}{n}}v; \sqrt{\frac{(k-1)!(n-1-k)!}{n}}v_k,\ k=\overline{1,n-1}\right\}$. Thus $$\begin{aligned}
		\det(\mathscr{C}_1(H))
		&=\det\left(G\left(\sqrt{\frac{\per(H)}{n}}v; \sqrt{\frac{(k-1)!(n-1-k)!}{n}}v_k,\ k=\overline{1,n-1}\right)\right)\\
		&=\frac{\per(H)}{n}\prod_{k=1}^{n-1}\frac{(k-1)!(n-1-k)!}{n}\cdot \det(G(v,v_1,\dotsc,v_{n-1}))\\
	\end{aligned}$$
And from the proof of proposition 3.4, we obtain that 
$$(v,v_1,\dotsc,v_{n-1})=(p_0,p_1,\dotsc,p_{n-1})\begin{pmatrix}
	1&\dotsc & ke_k&\dotsc &(n-1)e_{n-1}\\
	0&\dotsc& (-1)^2ne_{k-1}&\dotsc &(-1)^{2}ne_{n-2}\\
	.&.&.&.&.\\
	.&.&.&.&.\\
	.&.&.&.&.\\
	0&\dotsc&\dotsc&(-1)^jne_{k-j}&\dotsc\\
	.&.&.&.&.\\
	.&.&.&.&.\\
	.&.&.&.&.\\
	0&\dotsc&\dotsc&\dotsc&(-1)^{n-1}n\\
	
\end{pmatrix}$$
The matrix in the right side is the transition matrix given by
$$\begin{aligned}
	\text{The $(i,j)$-th entry}=
	\begin{cases}
		(-1)^ine_{j-i} & \text{if $i>1$}\\
		(j-1)e_{j-1} &\text{if $i=1$ and $j>1$}\\
		1 &\text{if $(i,j)=(1,1)$}\\

	\end{cases}
\end{aligned}$$
with convention that $e_0=1; \ e_t=0$ if $t<0$.
Moreover, we observe that the transition matrix is an upper triangular matrix with the absolute value of diagonal entries equal to $n$ except the $(1,1)$-th entry equal to 1 and $(p_0,p_1,\dotsc,p_{n-1})$ is a Vandermonde matrix. Hence 
$$\begin{aligned}
	\det(\mathscr{C}_1(H))
	&=\frac{\per(H)}{n}\prod_{k=1}^{n-1}n(k-1)!(n-1-k)!\cdot \det(G(p_0,p_1,\dotsc,p_{n-1}))\\
	&=\frac{\per(H)}{n}\prod_{k=1}^{n-1}n(k-1)!(n-1-k)!\cdot |\det(p_0,p_1,\dotsc,p_{n-1})|^2\\
	&=\frac{\per(H)}{n}\prod_{k=1}^{n-1}n(k-1)!(n-1-k)!\cdot \prod_{i<j}|x_i-x_j|^2\\
\end{aligned}$$
The right side is also equal to 0 if there are indices $i\neq j$ such that $x_i=x_j$. Hence the equality holds in both cases.
\end{proof}
\begin{remark}
	From the proposition 3.5, we are able to calculate the determinant of $\mathscr{C}_1(H)$ of any positive semi-definite Hermitian matrix $H$ of rank 2 in the way:\\
	Let $A$ be an $n\times n$ positive semi-definite Hermitian matrix of rank 2 then $A$ can be written in the form $vv^*+uu^*$ with $v_i,u_i$ are the $i$-th elements of $v$ and $u$ respectively. Then the following formula for the determinant of $\mathscr{C}_1(H)$ is achieved.\\
	\begin{theorem}
	Let $H=vv^*+uu^*$ be an $n\times n$ positive semi-definite Hermitian matrix then:
		$$\det(\mathscr{C}_1(H))=\frac{\per(H)}{n}\prod_{k=1}^{n-1}n(k-1)!(n-1-k)!\cdot \prod_{i<j}|v_iu_j-v_ju_i|^2$$
		where $v_i$ and $u_i$ are $i$-th elements of the vector $v$ and $u$ respectively.
		
	\end{theorem}

\end{remark}

\section{A counterexample for the conjectures 1 and 2 in the case $n=5$}
Let us take the values of $u_i$'s and $v_i$'s, $a\in \mathbb{R}$
$$u_1=ai,\ u_2=-a,\ u_3=-ai,\ u_4=a,\ u_5=0, \ v_i=1 \ \forall i=1,\dotsc,5 $$
then $e_1=e_2=e_3=e_5=0, \ e_4=-a^4$\\
For any matrix of the form, the spectrum of $\mathscr{C}_1(H)$ is determined clearly by the mentioned above properties and theorems.\\ 
By lemma 3.1, $\rank(\pi(H))\leq 2^5-5=27$ which means that there are at most 27 positive engenvalues.\\
By lemma 3.2, $$\per(H)=120+24|e_1|^2+12|e_2|^2+12|e_3|^2+24|e_4|^2+120|e_5|^2=120+24a^8$$
and the proposition 3.2 implies that $$\mathscr{C}_1(H)=\frac{\per(H)}{5}vv^*+\frac{6}{5}v_1v_1^*+\frac{2}{5}v_2v_2^*+\frac{2}{5}v_3v_3^*+\frac{6}{5}v_4v_4^*$$
where $$v_1=
\begin{pmatrix}
	-5ai\\
	5a\\
	5ai\\
	-5a\\
	0\\
\end{pmatrix}
\ v_2=
\begin{pmatrix}
	-5a^2\\
	5a^2\\
	-5a^2\\
	5a^2\\
	0\\
\end{pmatrix}
\ v_3=
\begin{pmatrix}
	5a^3i\\
	5a^3\\
	-5a^3i\\
	-5a^3\\
	0\\
\end{pmatrix}
\ v_4=
\begin{pmatrix}
	a^4\\
	a^4\\
	a^4\\
	a^4\\
	-4a^4\\
\end{pmatrix}$$
Notice that $\{v,v_1,v_2,v_3,v_4\}$ is orthogonal, thus those vectors are eigenvectors of $\mathscr{C}_1(H)$ corresponding to the eigenvalues $$\per(H)=120+24a^8,\frac{6}{5}\|v_1\|^2=120a^2, \ \frac{2}{5}\|v_2\|^2=40a^4,\ \frac{2}{5}\|v_3\|^2=40a^6,\ \frac{6}{5}\|v_4\|^2=24a^8$$
We replace $a^2=c,$ then $\tr(\pi(H))=120(1+c)^4$. The spectrum of  $\mathscr{C}_1(H)$ is $$\{120+24c^4, 120c,40c^2,40c^3,24c^4\}$$
 Moreover, every eigenvalue of  $\mathscr{C}_1(H)$ except $\per(H)$ is an eigenvalue of $\pi(H)$ with multiplicity at least 4 and, every eigenvalue of  $\mathscr{C}_2(H)$ except eigenvalues of $\mathscr{C}_1(H)$ is an eigenvalue of $\pi(H)$ with multiplicity at least 5. Thereforce, if we can calculate the sum and the sum of squares of at most 2 unknown positive eigenvalues of $\pi(H)$, then the spectrum is determined.
We compute the trace of $\mathscr{C}_2(H)$. The $(i,j)(i,j)$-th diagonal entry of $\mathscr{C}_2(H)$ is given by  $$\begin{aligned}
	& \ \per(H[\{i,j\},\{i,j\}]).\per(H(\{i,j\},\{i,j\}))\\
	& =(2+\left|e_1[\{i,j\}]\right|^2+2\left|e_2[\{i,j\}]\right|^2)(6+2\left|e_1[\{i,j\}^c])\right|^2+2\left|e_2[\{i,j\}^c]\right|^2+6\left|e_3[\{i,j\}^c]\right|^2)
\end{aligned}$$
Hence, we use the table to represent all the diagonal entries of $\mathscr{C}_2(H)$.
$$\begin{tabular}{|l|r|}
	\hline
	\text{Coordinates}&Values \\ \hline
	$(1,2)(1,2)$&$(2+2c+2c^2)(6+4c+2c^2)$ \\ \hline
	$(1,3)(1,3)$&$(2+2c^2)(6+2c^2)$ \\ \hline
	$(1,4)(1,4)$&$(2+2c+2c^2)(6+4c+2c^2)$ \\ \hline
	$(1,5)(1,5)$&$(2+c)(6+2c+2c^2+6c^3)$ \\ \hline
	$(2,3)(2,3)$&$(2+2c+2c^2)(6+4c+2c^2)$ \\ \hline
	$(2,4)(2,4)$&$(2+2c^2)(6+2c^2)$ \\ \hline
	$(2,5)(2,5)$&$(2+c)(6+2c+2c^2+6c^3)$ \\ \hline
	$(3,4)(3,4)$&$(2+2c+2c^2)(6+4c+2c^2)$ \\ \hline
	$(4,5)(4,5)$&$(2+c)(6+2c+2c^2+6c^3)$ \\ \hline
	$\tr(\mathscr{C}_2(H))$&$120+48c^4+104c^3+152c^2+120c$ \\ \hline
	
\end{tabular}$$
Furthermore, we use the symmetric polynomials to calculate the sum of all squares of eigenvalues.\\
$$\begin{aligned}
	\tr(\pi(H)^2)&=\sum_{\sigma \in S_5}\sum_{\tau \in S_5}\left|\prod_{i=1}^{5}(1+u_{\sigma(i)}\overline{u_{\tau(i)}})\right|^2\\
	&=120\sum_{\sigma \in S_5}\left|\prod_{i=1}^{5}(1+u_i\overline{u_{\sigma(i)}})\right|^2\\
\end{aligned}$$
We know that $u_5=0$, and for $k=1,\dotsc,4$ we have $u_k=a.i^k$ with $a^2=c$ then
$$\begin{aligned}
	\tr(\pi(H)^2)&=120\sum_{\sigma \in S_5}\left|\prod_{i=1}^{5}(1+u_i\overline{u_{\sigma(i)}})\right|^2\\
	&=120\left(\sum_{k=1}^{4}\sum_{\sigma \in S_5,\ \sigma(k)=5}\left|\prod_{j\neq k,5}(1+u_j\overline{u_{\sigma(j)}})\right|^2+\sum_{\sigma \in S_5,\ \sigma(5)=5}\left|\prod_{i=1}^{4}(1+u_i\overline{u_{\sigma(i)}})\right|^2\right)\\
	&=120\left(\sum_{k=1}^{4}\sum_{\sigma \in S_5,\ \sigma(k)=5}\left|\prod_{j\neq k,5}(1+c.i^{j-\sigma(j)})\right|^2+\sum_{\sigma \in S_5,\ \sigma(5)=5}\left|\prod_{j=1}^{4}(1+c.i^{j-\sigma(j)})\right|^2\right)\\
\end{aligned}$$
\begin{lemma}
	By the fundamental theorem of symmetric polynomials and $e_1=e_2=e_3=e_5=0$ then every monomial symmetric polynomial in 5 variables of degree non-divisible by 4 takes $(u_1,u_2,u_3,u_4,u_5)$ as a root.
\end{lemma}
The lemma 4.1 reduces the sums
$$\begin{aligned}
	&\quad \sum_{k=1}^{4}\sum_{\sigma \in S_5,\sigma(k)=5}\left|\prod_{j\neq k,5}(1+c\cdot i^{j-\sigma(j)})\right|^2\\
	&=\sum_{k=1}^{4}\sum_{\sigma \in S_5, \sigma(k)=5}(1+c^2)^3+(1+c^2)c\sum_{j \neq k,5}2\Rea(i^{j-\sigma(j)})\notag\\
	&\quad+(1+c^2)c^2\sum_{i_1<i_2 \neq k,5}(i^{i_1-\sigma(i_1)}+i^{\sigma(i_1)-i_1})(i^{i_2-\sigma(i_2)}+i^{\sigma(i_2)-i_2})+c^3\prod_{j\neq k,5}(i^{j-\sigma(j)}+i^{\sigma(j)-j})\\
	&=96(1+c^2)^3+\sum_{k=1}^{4}\sum_{\sigma \in S_5, \sigma(k)=5}c^2(1+c^2)2\Rea\left(\sum_{i_1<i_2 \neq k,5}i^{i_1-i_2+\sigma(i_2)-\sigma(i_1)}\right)\\
	&=96(1+c^2)^3+\sum_{k=1}^{4}c^2(1+c^2)\Rea\left(\sum_{i_1\neq i_2 \neq k,5}e^{i_1-i_2}\sum_{\sigma \in S_5,\ \sigma(k)=5}i^{\sigma(i_2)-\sigma(i_1)}\right)\\
\end{aligned}$$
combine with $$\sum_{\sigma \in S_5,\sigma(k)=5}i^{\sigma(i_2)-\sigma(i_1)}=2\sum_{\alpha=1}^{4}i^{\alpha}\sum_{\beta \neq \alpha}i^{\beta}=-2.4=-8$$\\
We attain $$\begin{aligned}
	\sum_{k=1}^{4}\sum_{\sigma \in S_5,\sigma(k)=5}\left|\prod_{j\neq k,5}(1+c\cdot i^{j-\sigma(j)})\right|^2&=96(1+c^2)^3-8c^2(1+c^2)\sum_{k=1}^{4}\Rea\left(\sum_{i_1\neq i_2 \neq k,5}e^{i_1-i_2}\right)\\
	&=96(1+c^2)^3+64c^2(1+c^2)\\
\end{aligned}$$
The lemma 4.1 also reduces the sum
$$\begin{aligned}
	&\sum_{\sigma \in S_5,\ \sigma(5)=5}\left|\prod_{i=1}^{4}(1+c\cdot i^{j-\sigma(j)})\right|^2=\sum_{\sigma \in S_4}\left|\prod_{i=1}^{4}(1+c\cdot i^{j-\sigma(j)})\right|^2\\
	&=\sum_{\sigma \in S_4}\left|1+c^4+c^3\sum_{i=1}^{4}i^{\sigma(j)-j}+c\sum_{i=1}^{4}i^{j-\sigma(j)}+c^2\sum_{j_1<j_2}i^{j_1+j_2-\sigma(j_1)-\sigma(j_2)}\right|^2\\
	&=24(1+c^4)^2+(c^6+c^2)\sum_{\sigma \in S_4}\left|\sum_{i=1}^{4}i^{j-\sigma(j)}\right|^2+c^4\sum_{\sigma \in S_4}\left|\sum_{j_1<j_2}i^{j_1+j_2-\sigma(j_1)-\sigma(j_2)}\right|^2\\
\end{aligned}$$
We compute each part separately by the lemma 4.1
$$\begin{aligned}
	&\sum_{\sigma \in S_4}\left|\sum_{i=1}^{4}i^{j-\sigma(j)}\right|^2=24\cdot4-8\sum_{j_1 \neq j_2}i^{j_1-j_2}=96+32=128\\
	&\sum_{\sigma \in S_4}\left|\sum_{j_1<j_2}i^{j_1+j_2-\sigma(j_1)-\sigma(j_2)}\right|^2\\
	&=\sum_{\sigma \in S_4}\biggl(\binom{4}{2}+\frac{1}{4}\sum_{\{i_1,i_2,i_3,i_4\}= \{1,2,3,4\}}i^{\sigma(i_1)+\sigma(i_2)-\sigma(i_3)-\sigma(i_4)+i_3+i_4-i_1-i_2}+2\sum_{j_1\neq j_2}i^{j_1-j_2+\sigma(j_2)-\sigma(j_1)}\biggr)\\
	&=144+2\sum_{(i_1,i_2,i_3,i_4)}i^{i_3+i_4-i_1-i_2}-16\sum_{j_1\neq j_2}i^{j_1-j_2}=208-4\sum_{j_1\neq j_2}i^{2j_1+2j_2}=224\\
\end{aligned}$$
 Thus, we obtain $\tr(\pi(H)^2)=120(24(1+c^4)^2+128(c^6+c^2)+224c^4+96(1+c^2)^3+64c^2(1+c^2))$.\\ Hence, the spectrum of $\pi(H)$ is
\begin{itemize}
	\item $\per(H)=120+24c^4$ of multiplicity 1
	\item $120c,40c^2,40c^3,24c^4$ of multiplicity 4
	\item $64c^3,112c^2$ of multiplicity 5
	\item 0 of multiplicity 93
\end{itemize} 
We observe that $c=2$ is a solution of the inequality $120+24c^4-64c^3<0$. Therefore, the matrix $H=vv^*+uu^*$ where $v=(1,\dotsc,1)^T,\ u=\sqrt{2}(i,-1,-i,1,0)^T$ is a counterexample to the permanent-on-top conjecture (POT).
	$$H=\begin{pmatrix}
	3&1-2i&-1&1+2i&1\\
	1+2i&3&1-2i&-1&1\\
	-1&1+2i&3&1-2i&1\\
	1-2i&-1&1+2i&3&1\\
	1&1&1&1&1\\
\end{pmatrix}$$
The spectrum of this counterexample is also given by above calculations:
	\begin{itemize}
	\item $\per(H)=504$ of multiplicity 1
	\item 240, 160, 320, 384 of multiplicity 4
	\item 512 and 448 of multiplicity 5
	\item 0 of multiplicity 93
\end{itemize}
Once, I have the counterexample, a shorter way to prove the matrix $H$ is a counterexample for Pate's conjecture in the case $n=5$ and $k=2$ is available by Tensor product. For the purposes of this paper let us describe the tensor product of vector spaces in terms of bases:
\begin{definition}
	Let $V$ and $W$ be vector spaces over $\mathbb{C}$ with bases $\{v_i\}$ and $\{w_i\}$, respectively. Then $V \otimes W$ is the vector space spanned by $\{v_i \otimes w_j\}$ subject to the rules:
	$$(\alpha v +\alpha ' v') \otimes w= \alpha (v \otimes w)+\alpha'(v' \otimes w)$$
	$$v\otimes (\alpha w+ \alpha ' w')=\alpha(v \otimes w)+ \alpha '(v \otimes w')$$
	for all $v, v' \in V$ and $w,w' \in W$ and all scalars $\alpha,\alpha '$.\\
	If $\langle,\rangle$ is an inner product on $V$ then we can define an inner product $\langle,\rangle$ on $V\otimes V$ in the manner: $$\langle v_{i_1}\otimes v_{i_2},v_{i_3}\otimes v_{i_4} \rangle=\langle v_{i_1},v_{i_3} \rangle \langle v_{i_2},v_{i_4} \rangle$$
	for any $v_{i_1},v_{i_2},v_{i_3},v_{i_4}$ vectors.
\end{definition}
On $\mathbb C[x,y]$, we consider the inner product, and the resulting Euclidean norm $|\cdot|$, such that monomials are orthogonal and $|x^ny^k|^2=n!k!$.
\begin{proposition}
	 The permanent of the Gram matrix of any 1-forms  $f_j\in \mathbb Cx\oplus \mathbb Cy$ is $\left|\prod f_j\right|^2$.
\end{proposition}
\begin{proof}
	We prove the generalization of the statement which states that if $f_1,f_2,\dotsc,f_n,g_1,g_2,\dotsc,g_n$ be $2n$ 1-forms and $A$ be an $n \times n$ matrix with $(i,j)$-th entry $\langle f_i,g_j \rangle$, then
	$$\per(A)=\left\langle \prod_{i=1}^{n}f_i, \prod_{i=1}^{n}g_i \right\rangle$$
	Let $f_i=\alpha_i x+ \beta_i y, g_i=\alpha'_i x+ \beta'_i y$ for any $i \in \{1,2,\dotsc,n\}$.\\
	We compute each side of the equality:\\
	The left side is
	$$\begin{aligned}
		\per(A)&=\sum_{\sigma \in S_n}\prod_{i=1}^{n}\langle f_i,g_{\sigma(i)} \rangle=\sum_{\sigma \in S_n}\prod_{i=1}^{n}\langle \alpha_i x+ \beta_i y,\alpha'_{\sigma(i)} x+ \beta'_{\sigma(i)} y \rangle=\sum_{\sigma \in S_n}\prod_{i=1}^{n}(\alpha_i \cdot \overline{\alpha_{\sigma(i)}}+\beta_i \cdot \overline{\beta'_{\sigma(i)}})\\
		&=\sum_{\sigma \in S_n}\sum_{k=0}^{n}\sum_{\substack{1\leq i_1<\dotsc<i_k\leq n \\ 1\leq i_{k+1}<\dotsc<i_n\leq n }}\alpha_{i_1}\dotsc \alpha_{i_k} \beta_{i_{k+1}}\dotsc\beta_{i_n}\overline{\alpha'_{\sigma(i_1)}\dotsc \alpha'_{\sigma{i_k}}}\overline{\beta'_{\sigma(i_{k+1})}\dotsc\beta'_{\sigma(i_n)}}\\
		&=\sum_{k=0}^{n}k!(n-k)!\left(\sum_{\substack{1\leq i_1<\dotsc<i_k\leq n \\ 1\leq i_{k+1}<\dotsc<i_n\leq n }}\alpha_{i_1}\dotsc \alpha_{i_k} \beta_{i_{k+1}}\dotsc\beta_{i_n}\right)\left(\sum_{\substack{1\leq i_1<\dotsc<i_k\leq n \\ 1\leq i_{k+1}<\dotsc<i_n\leq n }}\overline{\alpha'_1\dotsc \alpha'_{i_k} \beta'_{i_{k+1}}\dotsc\beta'_{i_n}}\right)\\
	\end{aligned}$$
and the right side is

	$$\begin{aligned}
		&\quad \left\langle \prod_{i=1}^{n}f_i, \prod_{i=1}^{n}g_i \right\rangle\\
		&=\left\langle\sum_{k=0}^{n}x^ky^{n-k}\sum_{\substack{1\leq i_1<\dotsc<i_k\leq n \\ 1\leq i_{k+1}<\dotsc<i_n\leq n }}\alpha_{i_1}\dotsc \alpha_{i_k} \beta_{i_{k+1}}\dotsc\beta_{i_n},\sum_{k=0}^{n}x^ky^{n-k}\sum_{\substack{1\leq i_1<\dotsc<i_k\leq n \\ 1\leq i_{k+1}<\dotsc<i_n\leq n }}\alpha'_{i_1}\dotsc \alpha'_{i_k} \beta'_{i_{k+1}}\dotsc\beta'_{i_n}\right\rangle\\
		&=\sum_{k=0}^{n}k!(n-k)!\left(\sum_{\substack{1\leq i_1<\dotsc<i_k\leq n \\ 1\leq i_{k+1}<\dotsc<i_n\leq n }}\alpha_{i_1}\dotsc \alpha_{i_k} \beta_{i_{k+1}}\dotsc\beta_{i_n}\right)\left(\sum_{\substack{1\leq i_1<\dotsc<i_k\leq n \\ 1\leq i_{k+1}<\dotsc<i_n\leq n }}\overline{\alpha'_1\dotsc \alpha'_{i_k} \beta'_{i_{k+1}}\dotsc\beta'_{i_n}}\right)\\		
	 \end{aligned}$$
\end{proof}
Let $f_j=x+yi^{j}\sqrt 2$  $(j=1,2,3,4)$ and $f_5=x$.  Their Gram matrix is the given matrix $H$ with $per H= |f_1f_2f_3f_4f_5|^2=|x^5-4xy^4|^2=5!+16\cdot 4!=504\ \text{(according to the proposition 4.1)}.$  
When   $\{p,q,r,s,t\}=\{1,2,3,4,5\}$, define $F_{p,q}=f_pf_q\otimes f_rf_sf_t$ and an inner product on $\mathbb{C}[x,y] \otimes \mathbb{C}[x,y]$ as the definition 4.1. It is obvious that $\mathscr{C}_2(H)$ of $H$ is the Gram matrix of the ten tensors $F_{p,q}$ with  $\{p,q,r,s,t\}=\{1,2,3,4,5\}$, $p<q$, and $r<s<t$. We observe that $$\begin{aligned}
	& \quad (1+i)F_{41}+(-1+i)F_{12}+(-1-i)F_{23}+(1-i)F_{34}-2iF_{51}+2F_{52}+2iF_{53}-2F_{54}\\ &=16\sqrt2x^2\otimes y^3-32\sqrt2 xy\otimes xy^2+16\sqrt2y^2\otimes x^2y,\end{aligned}$$
whose norm squared is $$2^9\cdot 2!3!+2^{11}\cdot 2!+2^9\cdot 2!\cdot 2!=512\cdot 24,$$ while the norm squared of the coefficient vector is $$|1+i|^2+|-1+i|^2+|-1-i|^2+|1-i|^2+|-2i|^2+2^2+|2i|^2+|-2|^2=24.$$
Therefore, a linear operator mapping eight orthonormal vectors to $F_{12}$, $F_{23}$, $F_{34}$, $F_{41}$, $F_{51}$, $F_{52}$, $F_{53}$, $F_{54}$  has norm at least $\sqrt{512}$, so the Gram matrix of these eight tensors, which is an 8-square diagonal submatrix
of $\mathscr{C}_2(H)$, has norm (=largest eigenvalue) at least 512, whence so does $\mathscr{C}_2(H)$  itself. In fact, the norm of $\mathscr{C}_2(H)$ is 512.\\

\end{document}